\documentclass[11pt]{article}
\usepackage{amscd,amssymb,stmaryrd}
\usepackage{amsmath}
\usepackage{amsthm}
\usepackage{color}
\usepackage{graphicx}
\usepackage{pifont}
\usepackage{subcaption}
\usepackage{tikz}
\usepackage[linesnumbered,ruled,vlined]{algorithm2e}

\usepackage{diagbox}
\usetikzlibrary{decorations.pathreplacing,calligraphy}
\usepackage{pgfplots}
\pgfplotsset{compat=1.11}
\usepgfplotslibrary{groupplots}

\usetikzlibrary{shapes.misc}

\tikzset{cross/.style={cross out, draw=black, minimum size=2*(#1-\pgflinewidth), inner sep=0pt, outer sep=0pt},
cross/.default={1pt}}
\theoremstyle{plain}
\theoremstyle{definition}

\newtheorem{remark}{Remark}

\def\ie{i.e., }

\newcommand{\kh}[1] {\textcolor{black}{#1}}
\title{BAYESIAN DEEP OPERATOR LEARNING FOR HOMOGENIZED TO FINE-SCALE MAPS FOR MULTISCALE PDE}
\author{ Zecheng Zhang\footnote{Department of Mathematics, Florida State University, Tallahassee, FL 32304, USA. (Email: zecheng.zhang.math@gmail.com)}, 
~Christian Moya \footnote{Department of Mathematics, Purdue University, West Lafayatte, IN 47907, USA. (Email: cmoyacal@purdue.edu )},
~Wing Tat Leung \footnote{Department of Mathematics, City University Hong Kong, Hong Kong, China. (Email: wtleung27@cityu.edu.hk)},
~Guang Lin \footnote{Department of Mathematics, Purdue University, West Lafayatte, IN 47907, USA. (Email: guanglin@purdue.edu  )},
~Hayden Schaeffer \footnote{Department of Mathematics, UCLA, Los Angeles, CA 90095. (Email: hayden@math.ucla.edu)} }
\begin{document}
\maketitle
\begin{abstract}
We present a new framework for computing fine-scale solutions of multiscale Partial Differential Equations (PDEs) using operator learning tools. Obtaining fine-scale solutions of multiscale PDEs can be challenging, but there are many inexpensive computational methods for obtaining coarse-scale solutions. Additionally, in many real-world applications, fine-scale solutions can only be observed at a limited number of locations. In order to obtain approximations or predictions of fine-scale solutions over general regions of interest, we propose to learn the operator mapping from coarse-scale solutions to fine-scale solutions using a limited number (and possibly noisy) observations of the fine-scale solutions. The approach is to train multi-fidelity homogenization maps using mathematically motivated neural operators.  The operator learning framework can efficiently obtain the solution of multiscale PDEs at any arbitrary point, making our proposed framework a mesh-free solver. We verify our results on multiple numerical examples showing that our approach is an efficient mesh-free solver for multiscale PDEs.
\end{abstract}

\section{Introduction}
Obtaining fine-scale solutions for a multiscale partial differential equation (PDE) problem can be costly, often requiring extensive large-scale computations to fully resolve.  In addition, while coarse scale solutions are easier to observe, we typically only have access to a limited number of fine-scale solutions in real-world multiscale applications. This limitation becomes particularly significant when fine-scale information is obtained by only a few samples or from specific subregions within the domain. Therefore, there is a need for methods that effectively utilize coarse scale solution (via simulations) alongside scarce fine-scale measurements. By addressing this one could compute accurate fine-scale solutions throughout the domain of interest.

Solving multiscale problems presents significant challenges due to the need for computationally expensive fine-scale solvers to accurately capture multiscale features. To address this issue, various multiscale methods have been proposed with the aim at resolving these complex phenomena. Among these methods, multiscale finite element methods (MsFEM) have emerged as particularly successful techniques \cite{efendiev2013generalized,chung2016adaptive, chung2018constraint, chung2022computational, chung2020generalized}. MsFEM operates by first solving for the multiscale basis, which effectively captures the local multiscale features on a coarse mesh. Subsequently, the problem is solved using this constructed basis, enabling a computationally effective approach for tackling multiscale problems. Although using MsFEM provides an effective framework for accurately representing and resolving the multiscale behavior inherent in these challenging problems, it can still be costly \textcolor{black}{\cite{chetverushkin2021computational, chung2022computational, leung2022nh, efendiev2022efficient, efendiev2023hybrid}}.

One technique for handling multiscale problems is through homogenization \cite{efendiev2009multiscale,allaire1992homogenization,leung2022nh}. A standard approach is to employ the asymptotic expansion to represent the solution. One can construct a homogenized solution without fully resolving the various scales, which still leads to a valid approximation to the exact solution. In particular, one can compute the solution to a homogenized PDE which only contains information at the coarse scale.  An alternative formulation is to establish a connection, in the form of an operator, between a solution at a finer scale and the homogenized solution.
The mathematical formulation and numerical computation of the operator for this task remains an open problem. However, we propose addressing this problem through operator learning by creating an efficient deep neural network-based approximation of the operator that represents the fine-scale solver of the corresponding multiscale PDE.

In \textit{operator learning}~\cite{lu2021learning,li2020fourier,zhang2022belnet, zhang2023discretization} one trains deep neural networks to approximate nonlinear operators, which are mappings between infinite-dimensional spaces. These operator-learning frameworks have been successful in scientific computing due to their versatility and efficiency for various problems arising from physical systems, including modeling dependencies on initial conditions or parameters. The first operator learning framework for PDE, the Deep Operator Network (DON), was developed in~\cite{lu2021learning}. DON is built on the universal approximation theorem for operators from ~\cite{chen1995universal} and can effectively learn operators with relatively small datasets. Compared to more traditional neural networks that learn mappings between vector spaces, DON exhibits improved generalization behavior on more complex tasks, as demonstrated in various applications. These include acting as surrogate solvers of PDEs (such as bubble dynamics~\cite{lin2021operator}), and approximating operators arising in tasks for control systems~\cite{li2022learning}, power grids~\cite{moya2023deeponet}, and multiphysics problems~\cite{cai2021deepm}. In \cite{zhang2022belnet, zhang2023discretization}, a discretization-invariant extension and analysis of DON was proposed, which allows the network to handle input functions with different discretizations. Some other extensions of DON have enabled incorporating physical information, leading to physics-informed DONs~\cite{wang2021learning}, handling noisy data~\cite{moya2023deeponet,lin2023b}, quantifying uncertainty~\cite{psaros2023uncertainty,moya2023deeponet,lin2023b, li2023fast}, or performing inverse design for complex applications~\cite{lu2022multifidelity, schaeffer2017learning}.

In a parallel effort, the Fourier Neural Operator (FNO) was presented in ~\cite{li2020fourier}. FNO approximates nonlinear operators by directly parametrizing integral kernels in the Fourier domain. FNO's effectiveness has been demonstrated in numerous applications across various domains, including but not limited to global weather prediction~\cite{pathak2022fourcastnet}, multiphase flow~\cite{wen2022u}, and solving PDEs with complex geometry~\cite{li2020fourier}. In~\cite{kovachki2021neural}, the authors generalized FNO and proposed neural operators that can effectively learn operators. Specifically, they formulated the neural operator as a composition of linear integral operators and nonlinear activation functions. Furthermore, the authors supported the neural operators by providing a universal approximation theorem, which demonstrates the existence of a neural operator that can approximate a given nonlinear continuous operator.

Training effective operator networks for multiscale problems can be expensive due to the requirement of large amounts of high-fidelity (fine-scale) data \cite{schaeffer2013sparse,schaeffer2017learning, lin2022theoretical}. However, in practice, high-fidelity data may be limited. Fortunately, we may have access to mathematical models that can generate abundant low-fidelity (coarse-scale) data. 
In such cases, our objective is to accurately calculate the fine-scale solution by constructing an operator learning framework that maps a function (i.e., the solution of the partially known PDE) to a new function that represents the exact solution of the target PDE. We then use the available data to train this framework, mixing scarce high-accuracy observation with abundant coarse simulations \cite{peherstorfer2018survey, peherstorfer2016optimal}.

Therefore, it is crucial to design operator learning frameworks that can be trained using high-fidelity data, mathematical models, and low-fidelity data. Several studies~\cite{howard2022multifidelity,lu2022multifidelity,de2023bi} have developed multi-fidelity neural operators to address this need, with applications to complex tasks such as fluid or materials science. However, they do not consider the uncertainty caused by using various fidelities or models during training and thus may not be effective in the presence of noisy data. To bridge this gap, this paper proposes a novel Bayesian deep operator learning architecture for developing surrogates of multi-scale PDE solvers. The architecture can incorporate noisy high-fidelity data and efficient solvers developed based on homogenization, which is a PDE-based technique for handling multiscale PDEs \cite{efendiev2009multiscale, allaire1992homogenization, robinson2008surrogate,hou2015sparse}. The main contributions of this paper are summarized below.
\begin{itemize}
    \item We propose a data-driven approach that downscales a given coarse model. This approach maps coarse-scale solutions to fine-scale solutions directly from data.
    \item We propose an ``oversampling'' approach to capture the input function to this operator through patches. In our numerical experiments, we observed that enlarging the patch leads to improvement in the prediction error.
    \item Furthermore, we have designed the first Bayesian, multiscale operator learning framework that is trained with noisy data. This framework can provide robust and mesh-free solutions, even from coarse-scale solutions.
    \item Finally, we demonstrate through multiple numerical experiments that the proposed framework represents an efficient mesh-free solver for multiscale PDEs.
\end{itemize}
The paper is organized as follows. Section~\ref{sec:problem-formulation} formulates the problem of learning the operator from coarse-scale solutions to fine-scale solutions. In Section~\ref{sec:DON-review}, we review the Deep Operator Network framework. Section~\ref{sec:proposed-methodology} provides detailed information about the proposed operator learning frameworks, which are trained with coarse-scale solutions and a limited number of observations of fine-scale solutions to approximate fine-scale solutions. Section~\ref{sec:Bayesian-DON} develops a Bayesian, multiscale operator learning framework that enables reliable prediction of fine-scale solutions, even when trained with noisy observations. We demonstrate the effectiveness of the proposed framework with a series of numerical examples in Section~\ref{sec:numerical-experiments}. Finally, Section~\ref{sec:conclusion} concludes the paper.

\section{Problem Formulation} \label{sec:problem-formulation}
The main goal is to improve the accuracy of a low-scale/low-accuracy model or physical simulation by using real observation data related to a specific physical multiscale process. To address this, the proposed framework involves first obtaining a coarse-scale solution with lower accuracy via a numerical solver. Then, an operator learning approach is used to refine the coarse-scale solution by incorporating available, possibly noisy, observed data. Obtaining the coarse-scale solution is often a more feasible task in terms of the computational cost or algorithmic complexity. Then, using the coarse-scale solution, the designed operator learning approach will act as a downscaled multiscale model, providing a fine-scale multiscale solver simultaneously.

We motivate and justify the descaling operator method as follows. Consider an example of homogenization of a multiscale elliptic problem:
\begin{align} \label{eq:homogeneization-multiscale}
    -\frac{\partial }{\partial x_i}\bigg( 
    a_{ij}(x/\epsilon) \frac{\partial}{\partial x_j} u_{\epsilon}(x)
    \bigg) = f(x), \quad x \in \Omega,
\end{align}
with $u_{\epsilon}(x) = 0$ at the boundary $\partial \Omega$. In the above equation, we have used the Einstein notation, {$u_\epsilon$ is the PDE solution}, {$a_{ij}$ is the multiscale permeability}, {$f(x)$ is the forcing}, and {$\Omega$ is the domain}. We seek $u_{\epsilon}(x)$ with an asymptotic expansion of the form:
\begin{align}
    u_{\epsilon}(x) = u_0(x, x/\epsilon)+\epsilon u_1(x, x/\epsilon)+
    \epsilon^2 u_2 (x, x/\epsilon) + \mathcal{O}(\epsilon^3)
    \label{eqn_homg_asym}
\end{align}
where $u_j(x, y)$, for $j=0,1,2,\ldots$, is periodic in $y = x/\epsilon$ with period $1$. Denote 
\begin{align*}
    A^\epsilon = -\frac{\partial }{\partial x_i}\bigg( 
    a_{ij}(x/\epsilon) \frac{\partial}{\partial x_j}\bigg).
\end{align*}
Then, it follows that:
\begin{align*}
    A^{\epsilon} = \epsilon^{-2} A_1+\epsilon ^{-1} A_2+ 
    \epsilon^0 A_3,
\end{align*}
where
\begin{align*}
    & A_1 = -\frac{\partial }{\partial y_i}\bigg( 
    a_{ij}(y) \frac{\partial}{\partial y_j}\bigg),\\
    & A_2 = -\frac{\partial }{\partial x_i}\bigg( 
    a_{ij}(y) \frac{\partial}{\partial y_j}\bigg)
    -\frac{\partial }{\partial y_i}\bigg( 
    a_{ij}(y) \frac{\partial}{\partial x_j}\bigg)\\
    & A_3= -\frac{\partial }{\partial x_i}\bigg( 
    a_{ij}(y) \frac{\partial}{\partial x_j}\bigg).
\end{align*}
Thus, we have the decomposition $\epsilon^{-2}A_1 u_{\epsilon} + \epsilon^{-1}A_2 u_{\epsilon} + \epsilon A_3 u_{\epsilon} = f$. By equating the terms with the same power of $\epsilon$, we obtain:
\begin{align}
    & A_1 u_0 = 0, \label{homog_u0} \\
    & A_1 u_1 + A_2 u_0 = 0, \label{homog_u1}\\
    & A_1u_2+A_2u_1+A_3u_0 = f. \label{homog_u2}
\end{align}
Substitute $A_1$ into equation \eqref{homog_u0} leads to:
\begin{align*}
    -\frac{\partial }{\partial y_i}\bigg( 
    a_{ij}(y) \frac{\partial}{\partial y_j}\bigg) u_0 = 0.
\end{align*}
According to the theory of second-order ordinary differential equations, $u_0$ is independent of $y$. This simplifies~\eqref{homog_u1} and we have:
\begin{align*}
    -\frac{\partial }{\partial y_i}\bigg( 
    a_{ij}(y) \frac{\partial}{\partial y_j}\bigg) u_1
    = \bigg( \frac{\partial }{\partial y_i}
    a_{ij}(y) \bigg) \frac{\partial}{\partial x_j} u_0.
\end{align*}
Thus, $u_1(x, y)$ can be solved by introducing $\chi_j(y)$, which is the solution to the following problem:
\begin{align*}
    &-\frac{\partial }{\partial y_i}\bigg( 
    a_{ij}(y) \frac{\partial}{\partial y_j}\bigg) \chi_j = 
    \frac{\partial}{\partial y_i}a_{ij}(y), \\
    & \chi_j \text{ is periodic in $y$ with mean $0$}.
\end{align*}
The equation above is referred to as a cell problem where it needs to be solved within one period of $y$, or, in the unit cell $Y = [0, 1]^d$, where $d$ is the dimension of the problem. Then, $u_1$ can be expressed as follows:
\begin{align*}
    u_1(x, y) = \chi_j \frac{\partial u_0}{\partial x_j}(x).
\end{align*}
Substituting the above into~\eqref{eqn_homg_asym}, we have:
\begin{align}
    u_{\epsilon, 1}(x) = u_0 + \epsilon\chi_j\frac{\partial u_0}{\partial x_j}(x),
\label{eqn_key_eqn}
\end{align}
which is a higher-order approximation of $u$ compared to $u_0$.

Finally, we can express \eqref{eqn_key_eqn} in operator form as follows:
\begin{align}
    u_{\epsilon, 1}(x) = G(u_0)(x).
\end{align}
Here the operator $G$ maps the {homogeneized solution} (\ie the \textit{coarse-scale solution}) $u_0$ to {a \textit{finer-scale solution}} $u_{\epsilon, 1}$. Our goal is to approximate the operator $G$ using a Deep Operator Network~(DON).
\section{A Brief Review of Deep Operator Networks} \label{sec:DON-review}
This section provides a review of the Deep Operator Network (DON) proposed in~\cite{lu2021learning}. DON is a neural network architecture that approximates mappings between infinite-dimensional spaces. It is built on the universal approximation theorem of continuous operators, which was introduced in the seminal works~\cite{chen1995universal, chen1993approximations}. In particular, DON satisfies the following approximation theorem. Suppose $X$ is a Banach space, $K_1\subset X$, and $K_2\subset \mathbb{R}$ are compact sets. If $V\subset C(K_1)$ is compact, then the continuous operator $G: V\rightarrow C(K_2)$ can be effectively approximated by a parameterized function. Specifically, for any $\epsilon>0$, there exist positive integers $M$, $N$, and $K$, constants $c_i^k$, $\zeta_k$, $\theta_i^k$, and $\varepsilon_{ij}^k\in\mathbb{R}$, and points $\omega_k\in\mathbb{R}^d$, $y_j\in K_1$, $i = 1, ..., M$, $k = 1, ..., K$, and $j = 1, ...., N$ such that 
\begin{align*}
        \bigg| G(u)(x) - \sum_{k = 1}^K \sum_{i = 1}^M c_i^k\, g\left(\sum_{j = 1}^N\varepsilon_{ij}^ku(y_j)+\theta_i^k\right)\, g(\omega_k\cdot x+\zeta_k)\bigg|<\epsilon
\end{align*}
holds for all $u\in V$ and $x\in K_2$.

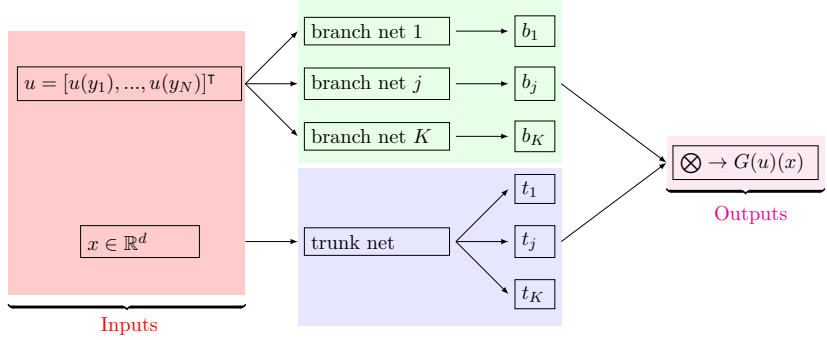
\begin{figure}[h!]
\centering
\scalebox{.7}{\begin{tikzpicture}[scale = 1]
 \fill [green!10] (2, 4.5) rectangle (7, 7.6);
  \fill [blue!10] (2, 1.4) rectangle (7, 4.4);
    \fill [red!20] (-3.5, 2) rectangle (1, 7);

    \draw[ultra thick] [decorate,
    decoration = {calligraphic brace, mirror}] (-3.5, 1.8) --  (1, 1.8);
\node at (-1.2, 1.4) {\textcolor{red}{Inputs}};

\fill [magenta!10] (9, 4) rectangle (12, 5);
    \draw[ultra thick] [decorate,
    decoration = {calligraphic brace, mirror}] (9, 4) --  (12, 4);
\node at (10.6, 3.5) {\textcolor{magenta}{Outputs}};

\node[draw, text width=4cm] at (-1.2, 6) {$u = [u(y_1), ..., u(y_N)]^\intercal$};

 \draw [-latex ](1,6) -- (2, 7);
 \node[draw, text width = 2.5cm] at (3.5, 7) {branch net $1$};
 \draw [-latex ](5, 7) -- (6, 7);
 \node[draw, text width = 0.5cm] at (6.5, 7) {$b_1$};
 
 \draw [-latex ](1,6) -- (2, 5);
  \node[draw, text width = 2.5cm] at (3.5, 5) {branch net $K$};
  \draw [-latex ](5, 5) -- (6, 5);
  \node[draw, text width = 0.5cm] at (6.5, 5) {$b_K$};
  
  \draw [-latex ](1,6) -- (2, 6);
 \node[draw, text width=2.5cm] at (3.5, 6) {branch net $j$};
 \draw [-latex ](5, 6) -- (6, 6);
\node[draw, text width = 0.5cm] at (6.5, 6) {$b_j$};

 \node[draw, text width = 2cm] at (-1, 3) {$x\in\mathbb{R}^d$};
  \draw [-latex ](1, 3) -- (2, 3);
 \node[draw, text width = 2.5cm] at (3.5, 3) {trunk net};
 
 \draw [-latex ](5, 3) -- (6, 4);
\node[draw, text width = 0.5cm] at (6.5, 4) {$t_1$};

  \draw [-latex ](5, 3) -- (6, 3);
\node[draw, text width = 0.5cm] at (6.5, 3) {$t_j$};

   \draw [-latex ](5, 3) -- (6, 2);
 \node[draw, text width = 0.5cm] at (6.5, 2) {$t_K$};

\draw [-latex ](7, 6) -- (9, 4.5);
  \draw [-latex ](7, 3) -- (9, 4.5);

  \node[draw, text width = 2.5cm] at (10.5, 4.5) {$\bigotimes\rightarrow G(u)(x)$};

\end{tikzpicture}}
\caption{Stacked version of the Deep Operator Network~(DON). $\bigotimes$ denotes the inner product in $\mathbb{R}^K$.}
\label{deepo_structure}
\end{figure}

The above approximation theorem suggests a neural network architecture for DON illustrated in Figure~\ref{deepo_structure}. The architecture comprises two sub-networks: a branch net and a trunk net. The \textit{branch net} is composed of a stacked collection of $K$ networks, which take the function $u$ discretized using $N$ sensors as input, and output the vector $(b_1,\ldots,b_K)^\top$. On the other hand, the \textit{trunk net} takes the location $x$ within the output function domain as input and outputs the vector $(t_1, \ldots, t_K)^\top$. The final output of the DON is obtained through the inner product between the output vectors of the branch and trunk nets.

It is important to note that the DON inputs contain the independent variable, $x$, which denotes the location of the output target function. This means that a well-trained DON can predict the output function value at any arbitrary point in its domain. We will use this property to construct an operator learning-based mesh-free solver for multiscale PDEs.
\section{Proposed Methodology} \label{sec:proposed-methodology}
We use DON to design two operator learning-based algorithms that can approximate the operator mapping coarse-scale solutions to fine-scale solutions.

\textit{The Vanilla Operator Learning Algorithm (without Patches)}. Here, we design and train a deep neural network $G_{\theta}$, with a vector of trainable parameters~$\theta$, to approximate the \textit{true} operator $G(u_0(x))(x)$, which maps the coarse-scale solution $u_0(x)$ to the fine-scale solution at any given location~$x$ within the domain~$\Omega$. We summarize this vanilla algorithm in Algorithm \ref{algo_vanilla}. We also note that we refer to this algorithm as ``without patch'' to emphasize that we only use the coarse solution at a point $x$, and not any neighboring locations within $\Omega$.

To train the proposed DON $G_\theta: u_0 \mapsto u$, we minimize the loss function:
$$\mathcal{L}(\theta) = \frac{1}{N_p} \sum_{i=1}^{N_p} \|G_\theta(u_0(x_i))(x_i) - u(x_i)\|^2,$$
using the dataset of $N_p$ triplets $\{u_0(x_i), x_i, u(x_i)\}_{i=1}^{N_p}$, where $u_0(x_i)$ is the low-accuracy (coarse-scale) solution, $x_i \in \Omega$ the given location within the domain~$\Omega$, and $u(x_i)$ the fine-scale (observation) solution.

\begin{algorithm}[h!]
\caption{Vanilla Operator Learning Without Patch} \label{algo_vanilla}
    \textbf{Requires:} dataset: $\{u_0(x_i), x_i, u(x_i)\}_{i=1}^{N_p}$, where $u_0(x_i)$ is the low-accuracy (coarse-scale) solution, $x_i$ the location within the domain~$\Omega$, and $u(x_i)$ the fine-scale (observation) solution. \\
    Use a coarse solver to solve the multiscale equation and obtain $u_0(x_i)$, where $x_i \in \Omega$.\\
    Train the DON $G_{\theta}: u_0 \mapsto u$ by minimizing the loss function
\begin{align} \label{eqn_loss}
    \mathcal{L}(\theta) = \frac{1}{N_p}\sum_{i = 1}^{N_p}\|G_{\theta}(u_0(x_i))(x_i) - u(x_i)\|^2.
\end{align}\\
    Predict fine-scale solutions at arbitrary locations $x\in\Omega$ using the trained DON $G_{\theta^*}(u_0)(x)$, where $\theta^*$ denotes the optimal network parameters.
\end{algorithm}

\begin{remark}
    We have two remarks regarding Algorithm~\ref{algo_vanilla}.
    \begin{enumerate}
        \item The coarse-scale solution $u_0$ can be obtained easily and with low computational cost. One approach is to use multiscale finite element methods (MsFEM) \cite{efendiev2013generalized, chung2016adaptive, chetverushkin2021computational}. Another option is to use data-free machine learning methods, such as neural homogenization physics-informed neural networks (NH-PINN) \cite{leung2022nh}. PINN is a mesh-free solver, which means it can calculate the solution at any point in the PDE's domain, allowing $u_0$ to be evaluated at any point in the domain.
        \item The resulting trained DON, $G_{\theta^*}(u_0(x))(x)$, with optimal parameters $\theta^*$, is a fine-scale, mesh-free solver.
    \end{enumerate}
\end{remark}

\textit{The Operator Learning Algorithm with patch}. In the vanilla algorithm, the input function $u_0$ (\ie the coarse-scale solution) is evaluated at a single point $x_i$. However, in {Equation \ref{eqn_key_eqn}}, $\frac{\partial u_0}{\partial x}(x_i)$ requires the consideration of the derivative at that point. To approximate this derivative, finite difference schemes use neighboring points, leading to the idea of sampling a local neighborhood centered around $x_i$, \ie \textit{the patch}.

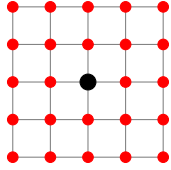
\begin{figure}[h!]
    \centering
\begin{tikzpicture}[scale=1]

\draw[step=0.5cm, gray] (-4,0) grid (-2., 2);

\foreach \i in {0, ..., 4}
{
\foreach \j in {0, ..., 4}
{
\filldraw[red] (-4+\i/2, \j/2) circle (2pt) node[anchor = west]{};
}
}

\filldraw[black] (-3, 1) circle (3pt) node[anchor=west]{};

\end{tikzpicture}
    \caption{A $5\times 5$ patch (red dots) centered at an observation point (black dot). All red dots represent the sensors used to sample the local input function (coarse-scale solutions) centered at the black dot (observation).}
    \label{fig_patch_domon}
\end{figure}

\begin{figure}[t!]
\centering
\scalebox{.71}{\input{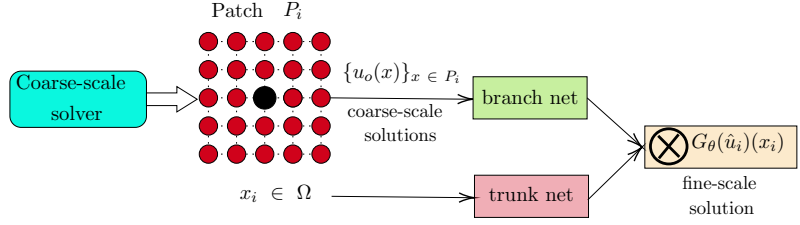}}
\caption{
A pictorial description of the operator learning algorithm with the patch. First, a coarse solver generates a collection of coarse-scale solutions $\hat{u}_i = \{u_o(x)\}_{x \in P_i}$ on the patch $P_i$ centered at $x_i \in \Omega$. The collection of coarse-scale solutions $\hat{u}_i$ is then input to the branch net, while the target location $x_i \in \Omega$ is input to the trunk net. The proposed operator learning algorithm with patch outputs the fine-scale solution $G_\theta(\hat{u}_i)(x_i) = u(x_i)$ at the arbitrary target location $x_i \in \Omega$.
}
\label{fig:patch-algorithm}
\end{figure}

To illustrate this concept, we provide an example (see Figure \ref{fig_patch_domon}) involving a $5\times 5$ patch for a two-dimensional (2D) problem. Subsequently, we propose a novel algorithm, the operator learning with patch algorithm detailed in Algorithm~\ref{algo_patch} and illustrated in Figure~\ref{fig:patch-algorithm}, which uses the patch sampling approach to estimate the desired derivative. 

Formally, we approximate the fine-scale solution $u(x_i)$ at any arbitrary point $x_i$ within the domain using the DON $G_\theta(\hat{u_i})(x_i)$. Here, $\hat{u}_i$ refers to the collection of coarse-scale solutions $\{u(x)\}_{x \in P_i}$ evaluated at the patch $P_i$ (i.e., the neighborhood of locations around $x_i$), which serves as the new input to the branch net.

To train the proposed DON $G_\theta: \hat{u}_i \mapsto u(x_i)$, we minimize the loss function:
$$\mathcal{L}(\theta) = \frac{1}{N_p} \sum_{i=1}^{N_p} \|G_\theta(\hat{u}_i)(x_i) - u(x_i)\|^2,$$
using the dataset of $N_p$ triplets $\{\hat{u}_i, x_i, u(x_i)\}_{i=1}^{N_p}$. 

\begin{algorithm}[h!]
\caption{Operator Learning With Patch} \label{algo_patch}
    \textbf{Require}: dataset: $\{\hat{u}_i, x_i, u(x_i)\}_{i=1}^{N_p}$. Here $u(x_i)$ is the fine-scale (observation) solution and 
    $\hat{u}_i = \{u(x)\}_{x\in P_i}$ is  the collection of coarse-scale solutions $\{u(x)\}_{x \in P_i}$ evaluated at the patch $P_i$ (i.e., the neighborhood of locations around $x_i$). For example, a patch of size three around $x_i$ of 1D case is $P_i = \{x_{i-1}, x_i, x_{i+1}\}$. A $5\times 5$ 2D patch demonstration is presented in Figure \ref{fig_patch_domon}.\\
    Use a coarse solver to solve the multiscale equation and obtain $u_0(x)$, where $x \in P_i$. \\
    Train the DON $G_{\theta}: \hat{u}_i \mapsto u(x)$ by minimizing the loss function
\begin{align*}
    \mathcal{L}(\theta) = \frac{1}{N_p}\sum_{i = 1}^{N_p}\|G_{\theta}(\hat{u}_i)(x_i) - u(x_i)\|^2.
\end{align*}\\
    Predict fine-scale solutions at $x_i\in\Omega$ using the trained DON $G_{\theta^*}(\hat{u}_i)(x_i)$, where $\theta^*$ denotes the optimal network parameters.
\end{algorithm}

\begin{remark}
    We conclude this section with two remarks concerning the patch.
    \begin{enumerate}
        \item Our numerical experiments have revealed a notable trend: as the size of the patch or neighborhood increases, the relative error decreases.
        \item We use DON to learn the operator. However, since DON is not invariant to the discretization of the input function, all patches must share the same discretization.
    \end{enumerate}
\end{remark}
\section{Bayesian, Multiscale Operator Learning} \label{sec:Bayesian-DON}
The conventional optimization framework used to train DON does not accurately quantify uncertainty and provide robust predictions, which are crucial for creating credible intervals for scientific and engineering applications. This can lead to unreliable DON predictions and limit the reliability of DON. However, quantifying the uncertainty associated with limited and noisy training data, along with the fact that the network is over-parametrization makes this a challenging task. And this challenge is even greater in operator learning, as it involves mappings between infinite-dimensional spaces, such as the multiscale operator~$G$. In this section, we address this challenge by developing a Bayesian DON (B-DON)~\cite{lin2023b}. B-DON can create estimators and credible intervals for the operator that maps the homogenized solution $u_o$ to the fine-scale solution $G(u_o)(x)$ for any given $x \in \Omega$.

In this Bayesian DON framework, our goal, given the training noisy dataset $\mathcal{D}$, is to construct a distribution $p(G|(u_o,x),\mathcal{D})$ that can predict the operator value of $G$ (the fine-scale solution) based on the input homogenized solution $u_o$ and at any new location $x$. To achieve this, we first assume the following factorized Gaussian likelihood function for the data:
\begin{align} \label{eq:gaussian-likelihood}
p(G|(u_o,x),\theta) = \mathcal{N}(G|G_\theta(u_o)(x), \text{diag}(\Sigma^2)) = \prod_{j=1}^N \mathcal{N}(G_j|G_\theta(u_{o,j})(x_j), \sigma),
\end{align}
where the output $G_\theta(u_o)(x)$ is the mean of the Gaussian distribution assumed for $G$, given the homogenized input $u_o$ at location $x$, and $\text{diag}(\Sigma^2)$ is a diagonal covariance matrix with $\Sigma^2 = (\sigma^2,\ldots,\sigma^2)$ on the diagonal~\cite{psaros2023uncertainty}. Note that~$\sigma$ can be assumed or estimated from the noisy data.

The fine-scale solution $G$ for a given homogeneized solution $u_o$ at a location $x$, given the noisy training data $\mathcal{D}$, is the random variable $(G|(u_o,x),\mathcal{D})$. To obtain the density of this random variable, we need to integrate the model parameters as follows:
$$p(G|(u_o,x), \mathcal{D}) = \int p(G|(u_o,x),\theta)p(\theta | \mathcal{D})d\theta.$$
Here, $p(\theta | \mathcal{D})$ represents the posterior distribution of the trainable parameters. This distribution enables us to quantify the \textit{epistemic uncertainty}, which refers to the uncertainty related to the trainable parameters $\theta$~\cite{psaros2023uncertainty,lin2023b}.

To obtain this posterior, we use Bayes' rule:
$$p(\theta|\mathcal{D}) \propto p(\mathcal{D}|\theta) p(\theta),$$
where $p(\theta)$ is the \textit{prior} distribution of the parameters and $p(\mathcal{D}|\theta)$ the \textit{data likelihood}, i.e., $p(\mathcal{D}|\theta) = \prod_{j=1}^N p(G_j|(u_{o,j}, x_j), \theta)$, which we calculate using the DON forward pass and the i.i.d noisy training dataset~$\mathcal{D}$. Acquiring the posterior distribution using Bayes' rule is computationally and analytically intractable~\cite{psaros2023uncertainty}. Therefore, in previous work~\cite{lin2023b}, we approximated the posterior distribution through samples obtained from it. Specifically, we obtained an $M$-ensemble of $\theta$ samples, denoted as $\{\theta_k\}_{k=1}^M$, as described below.
\subsection{Sampling the $M$-ensemble $\{\theta_k\}_{k=1}^M$}  \label{subsec:sampling-M-ensemble}
To obtain the $M$-ensemble of parameters $\{\theta_k\}_{k=1}^M$, B-DON uses the stochastic gradient replica exchange Langevin diffusion (SG-reLD), which we developed and studied in~\cite{lin2023b,lin2022multi,lin2021accelerated,na2022replica,li2023fast, deng2020non}. As demonstrated in our previous work, SG-reLD enjoys theoretical guarantees beyond convex scenarios, effectively handles large datasets, and accelerates convergence to the posterior distribution $p(\theta|\mathcal{D})$.

Specifically, SG-reLD uses two Langevin diffusions to describe the stochastic dynamics of $\theta$, along with a stochastic process that allows the diffusions to swap simultaneously. The high-temperature diffusion enables exploration of the parameter space, facilitating convergence to the flattened distribution of $\theta$. The low-temperature diffusion exploits the same parameter space, enabling faster convergence to local minima $\theta^*$. By swapping the diffusions, SG-reLD effectively escapes local minima and allows $\theta_k$ to converge faster to the desired posterior $p(\theta|\mathcal{D})$. For more details about employing SG-reLD with B-DONs, please refer to our previous paper \cite{lin2023b,lin2022multi}.

In practice, we can use the $M$-ensemble $\{\theta_k\}_{k=1}^M$ obtained using SG-reLD to fit a parametric distribution, such as the Gaussian distribution $$\mathcal{N}(\bar{\mu}(u_o)(X_{\text{test}}), \bar{\sigma}^2_e(u_o)(X_{\text{test}}))$$ for an arbitrary mesh~$X_{\text{test}}$ of locations. The parameters of this distribution are given by:
\begin{align}
    \bar{\mu}(u_o)(X_{\text{test}}) &= \frac{1}{M} \sum_{k=1}^M G_{\theta_k}(u_o)(X_{\text{test}}), \label{eq:mean} \\
    \bar{\sigma}^2_e(u_o)(X_{\text{test}}) &= \frac{1}{M}\sum_{k=1}^M \left(G_{\theta_k}(u_o)(X_{\text{test}}) - \bar{\mu}(u_o)(X_{\text{test}})  \right)^2. \label{eq:std}
\end{align}
Sampling from the aforementioned distribution allows for estimating credible sets for the fine solution. This also enables the sampling of more reliable predictions, which can help reduce the average relative error of test trajectories in the presence of noise, as demonstrated in our numerical experiments section.
\section{Numerical experiments}
\label{sec:numerical-experiments}
This section presents several numerical experiments to demonstrate the effectiveness of the proposed framework. For all examples, we conducted 100 independent experiments and presented the average relative error. We will analyze the error decay with respect to the patch size, the error decay with respect to the number of observation points, and the error decay in the presence of noise.
\subsection{1D Elliptic}
\begin{figure}[h!]
\centering
\includegraphics[scale = 0.5]{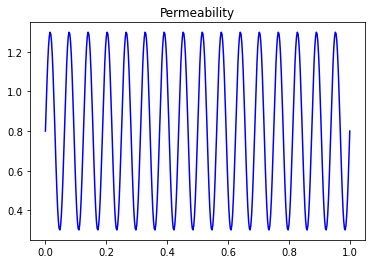}
\includegraphics[scale = 0.5]{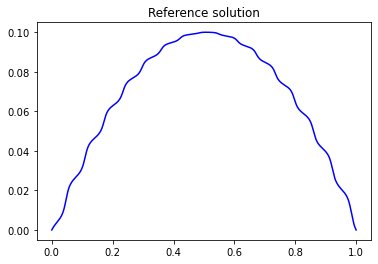}
\caption{1D elliptic. \textit{Left:} permeability $\kappa (x)$. \textit{Right:} reference solution.}
\label{1d_sprb}
\end{figure}

\begin{figure}[h!]
\centering
\includegraphics[scale = 0.5]{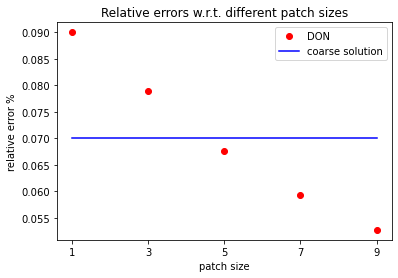}
\caption{Relative errors for 1D elliptic equations are shown with respect to different patch sizes. The patch size ranges from $1$ (using only the observation point) to $9$ (using $9$ points with the observation point in the patch center). We trained $100$ independent models and present their average relative errors. }
\label{fig_1d_results}
\end{figure}

In our first example, we will study a 1D problem for which we can obtain an exact homogenized solution. Specifically, we will consider the following elliptic equation:
\begin{align*}
    &-\frac{d}{dx} (a(x/\epsilon)\frac{du}{dx}) = f, \quad x\in[0, 1],\\
    & u(0) = u(1) = 0,
\end{align*}
where $a(x) = 0.5\sin(2\pi\frac{x}{\epsilon}) + 0.8$ and $f(x) = 0.5$. Figure~\ref{1d_sprb} illustrates the multiscale permeability $\kappa(x)$ and the reference solution.
We set $\epsilon = 1/16$. Then, the coarse-scale solution is obtained using classical homogenization theory, and the relative error of the homogenized solution is $0.07\%$. To further improve the homogenized solution, we employ data from exact solutions at uniformly distributed $16$ points in the domain, i.e., $N_p = 16$ in~\eqref{eqn_loss}. We evaluate the performance of our approach as the patch enlarges using the oversampling trick, and we present the results in Figure~\ref{fig_1d_results}.
\subsection{2D Elliptic Equation with $1$ Fast Variable}
\begin{figure}[h!]
\centering
\includegraphics[scale = 0.5]{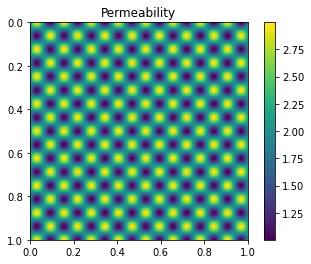}
\includegraphics[scale = 0.5]{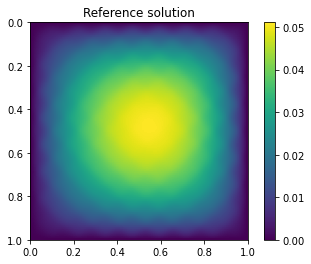}
\caption{2D elliptic with $1$ fast variable. \textit{Left:} permeability $\kappa(x)$. \textit{Right:} reference solution.}
\label{2d_kappa}
\end{figure}

In this experiment, we consider the following 2D elliptic equation:
\begin{align}
    &-\nabla\cdot (\kappa(x/\epsilon)\nabla u) = f, x\in\Omega = [0, 1]^2,\\
    & u(x) = 0, x\in \partial\Omega,
    \label{eqn_elliptic}
\end{align}
where $\kappa(x/\epsilon) = 2 + \sin(2\pi x/\epsilon)\cos(2\pi y/\epsilon)$ and $\epsilon = \frac{1}{8}$. Figure~\ref{2d_kappa} shows the permeability $\kappa$ as a function of $x$.
For the proposed framework, we assume prior knowledge of the exact solution of the given equation at a limited number of points. To compute the coarse-scale solution, we use the methodology described in \cite{leung2022nh}. Specifically, we use a neurohomogenized physics-informed neural network (NH-PINN) to derive the homogenized solution, denoted as $u_0$. The NH-PINN method employed in this study is a mesh-free solver, which allows for the generation of coarse-scale solutions at any point within the domain. As a result, we can always obtain a patch consisting of coarse-scale solutions centered around the observed fine-scale solutions.

Our goal for this experiment is to use operator learning to approximate the mapping from the NH-PINN-based coarse-scale solution $u_0$ to the corresponding fine-scale solution $u_f$. By learning this mapping, we aim to improve the accuracy of the coarse-scale solution by leveraging the information from the fine-scale solution.

\begin{figure}[h!]
\centering
\includegraphics[scale = 0.5]{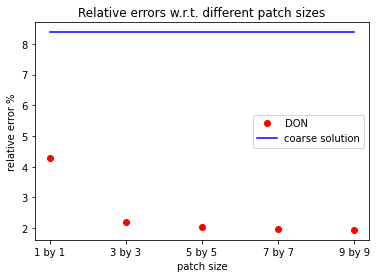}
\caption{We investigate the relative errors associated with different patch sizes, ranging from $1 \times 1$ (using only the observation point) to $9 \times 9$ (consisting of a total of $81$ points, with the observation point located at the center of the patch). To obtain a thorough evaluation, we train a total of $100$ independent models and report the average relative errors.}

\label{fig_2d_results}
\end{figure}

We demonstrate that the operator can be constructed (trained) more effectively as the patch size increases. We conducted five sets of experiments with patch sizes of $1\times 1$, $3\times 3$, $5\times 5$, $7\times 7$, and $9\times 9$. For each set of experiments, we trained $100$ models and computed the average relative errors of the last $100$ epochs for each model. The results are shown in Figure~\ref{fig_2d_results}.

Remarkably, our analysis shows that as the patch size increases, there is a reduction in relative errors. This observation highlights the benefits of enlarging the patch size in order to achieve improved accuracy in our models.

In the second part of this numerical experiment, we investigate the relationship between the relative error and the number of observation (training) points, which represent the exact solution. To achieve this, we use a fixed patch size of $1$, which means that we only include one point from the neighborhood for each observation sample.

\begin{figure}[h!]
\centering
\includegraphics[scale = 0.5]{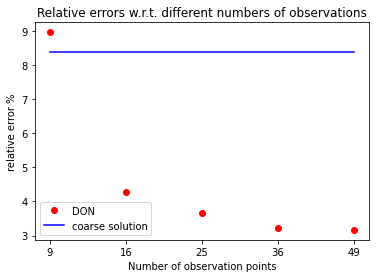}
\caption{We use various exact solutions, uniformly distributed in the domain, as the training labels (observations). We set the patch size to be $1$, which means that only the observation point coordinate is included as the branch input. By increasing the number of observations from $9$ to $49$, we evaluate the predictions on a $100\times 100$ mesh to compare the errors. For each observation set, we train a total of $100$ independent models, which capture the inherent variability in the training process. To assess the performance, we calculate the average relative error across the ensemble of models. This approach provides a comprehensive evaluation of the predictive accuracy for different numbers of observations.}
\label{fig_2d_wrt_obs}
\end{figure}

To investigate the influence of the number of observation points, we vary the number of observation points and present the results in Figure \ref{fig_2d_wrt_obs}. This analysis helps us understand the impact of the number of observation points on the relative error, providing valuable insights into the behavior of our model.

\begin{figure}[h!]
\centering
\includegraphics[scale = 0.5]{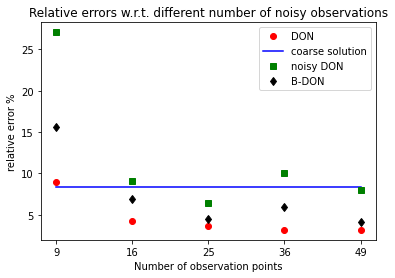}
\caption{We employ different numbers of exact solutions as the training labels (observations). We add Gaussian noise with variance $\sigma^2 = 0.005^2$ to the training labels. We fix the patch size to be $1$, i.e., include only the observation point coordinate as the brunch input.
The number of observations runs from $9$ to $49$, while we test the prediction on $100\times 100$ mesh. For each set of observations, we use as reference the results from the previous section that trains $100$ models with the true labels (DON). We train another $100$ models with noisy targets (noisy DON). We also use the proposed Bayesian framework (B-DON) to sample from the predictive distribution constructed using an ensemble of $M=100$ sets of parameters. For all the cases, we compute the average relative error.}
\label{fig_2d_wrt_noisy_obs}
\end{figure}

To conclude this experiment, we train the proposed multiscale B-DON model using noisy observations. We again conduct experiments by gradually increasing the number of observations from 9 to 49, while constructing an M-ensemble of models. Then, by sampling the fitted distribution whose parameters are given in~\eqref{eq:mean} and~\eqref{eq:std} and using a $100 \times 100$ $X_\text{test}$ mesh, we predict the values of the fine-scale solutions. Figure~\ref{fig_2d_wrt_noisy_obs} depicts the results of such experiments, which illustrate that the proposed Bayesian multi-fidelity operator learning framework can provide robust predictions even in the presence of noisy observations.

\subsection{2D Elliptic with Multiple Scales}
\begin{figure}[h!]
\centering
\includegraphics[scale = 0.5]{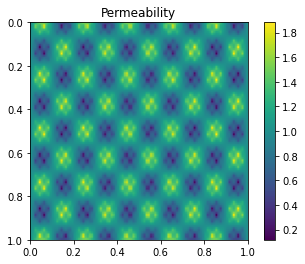}
\includegraphics[scale = 0.5]{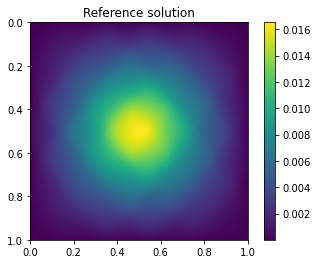}
\caption{2D elliptic with multiple scales. \textit{Left:} Permeability $\kappa$. \textit{Right:} Reference solution.}
\label{non_sprb}
\end{figure}

This experiment considers the same equation~\eqref{eqn_elliptic} as before but with a different permeability $\kappa$. Specifically, we let $\kappa$ be:
\begin{align*}
    \kappa(x, y) = 1 + \frac{\sin(2\pi\frac{x}{\epsilon_0})\cos(2\pi\frac{y}{\epsilon_1}) }{2+\cos(2\pi\frac{x}{\epsilon_2})\sin(2\pi\frac{y}{\epsilon_3})} +
    \frac{\sin(2\pi\frac{x}{\epsilon_4})\cos(2\pi\frac{y}{\epsilon_5})}{2+\cos(2\pi\frac{x}{\epsilon_6})\sin(2\pi\frac{y}{\epsilon_7})},
\end{align*}
where $\epsilon_0 = \frac{1}{5}$, $\epsilon_1 = \frac{1}{4}$, $\epsilon_2 = \frac{1}{25}$, $\epsilon_3 = \frac{1}{16}$, $\epsilon_4 = \frac{1}{16}$, $\epsilon_5 = \frac{1}{32}$, $\epsilon_6 = \frac{1}{3}$, $\epsilon_7 = \frac{1}{9}$. Figure~\ref{non_sprb} illustrates the permeability and reference solution.

\begin{figure}[h!]
\centering
\includegraphics[scale = 0.5]{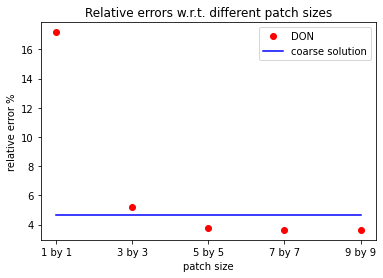}
\caption{Relative errors of the 2D elliptic example with respect to different patch sizes. The patch size ranges from $1\times1$ (using only the observation point) to $9\times9$ (using $81$ points with the observation point at the patch center). We trained $100$ independent models and present the average relative errors.}
\label{fig_2d_results2}
\end{figure}

We obtain the coarse-scale solution using multiscale finite element methods with one local basis \cite{efendiev2013generalized,chung2016adaptive,chung2018constraint,chetverushkin2021computational}. Then, we demonstrate that the approximation to the true solution operator can be better constructed (learned) as the patch size increases. To illustrate this, we conduct three sets of experiments with patch sizes of $1\times 1$, $3\times 3$, $5\times 5$, $9\times 9$, and $16\times 16$. For each set of experiments, we trained $100$ models and computed the average relative errors of the last $100$ epochs of all the $100$ models. Figure \ref{fig_2d_results2} shows the obtained results.
\subsection{Radiative Transfer Equation}
\begin{figure}[h!]
    \centering
    \includegraphics[scale = 0.5]{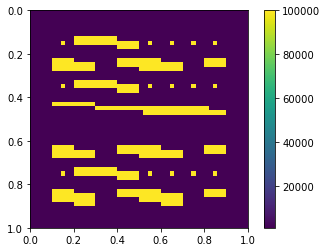}
    \includegraphics[scale = 0.5]{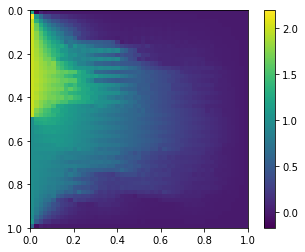}
    \caption{\textit{Left:} demonstration of multiscale scattering $\sigma(x)/\epsilon$, where $\epsilon = 0.001$. \textit{Right:} the solution of the RTE.}
    \label{exp_demon_sampling}
\end{figure}

\begin{figure}[h!]
    \centering
    \includegraphics[scale = 0.5]{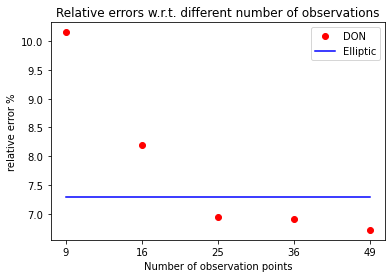}
    \caption{We use varying numbers of RTE solutions, uniformly distributed in the domain, as training labels (observations). The number of observations ranges from $9$ to $49$. We test the prediction on a $51 \times 51$ mesh. For each set of observations, we train $50$ independent models and calculate the average relative error.}
    \label{demon_rte_results}
\end{figure}

In this final experiment, we consider the Radiative Transfer Equation (RTE) \cite{newton2020diffusive, li2019diffusion, chung2020generalized} with a high-contrast scattering coefficient $\sigma(x, \omega)$ (see Figure \ref{exp_demon_sampling}). The term `high-contrast' refers to the strong scattering in the channels:
\begin{align*}
 s \cdot \nabla  I(x,s) = \frac{\sigma(x, \omega)}{\epsilon} \bigg( \int_{\mathcal{S}^{n-1}} I(x,s') ds' - I(x, s)\bigg) \quad \forall x \in D, s \in \mathcal{S}^{n-1}.
\end{align*}
Here, $s$ is a vector on the unit sphere, and $n$ is the dimension of the problem. In our experiments, we considered $n=2$, and thus $\mathcal{S}^{n-1}=\mathcal{S}^{1}$ is the unit circle. Additionally, we set $\epsilon = 0.001$ and $D = [0, 1]^2$. We also introduced the Dirichlet boundary conditions $I(x,s) = I_{\text{in}}$ for entrant directions $s \cdot \textbf{n} <0$, i.e., on $\Gamma ^{-} := \{ (x,s) \in \partial D \times \mathcal{S}^{n-1}: s \cdot \textbf{n} <0 \}$. Here, $\textbf{n}$ is the unit outward normal vector field at $x \in \partial D$. The condition can be written as:
\begin{align*}
I = I_{\text{in}} (x,s) \quad \text{ for all } (x,s) \in \Gamma^-.
\end{align*}
In our examples, the top, bottom, and right boundaries have zero incoming boundary conditions. We also assume that the left boundary has non-zero flow injected into the domain.
We chose the multiscale Radiative Transfer Equation (RTE) with high contrast channels for its numerical complexity \cite{chung2020generalized} and its challenge for learning-based approaches. However, as $\epsilon$ approaches zero, the elliptic solution converges to the RTE \cite{newton2020diffusive, li2019diffusion}. We use the elliptic solution as a low-accuracy solution and incorporate observed real RTE solutions to learn the downscaling of the model. Specifically, the RTE solution is used to correct errors in the elliptic solution. To understand how observations improve downscaling, we conducted a series of experiments with different numbers of observation points. We present the results in Figure~\ref{demon_rte_results}.
\section{Conclusion} \label{sec:conclusion}
\kh{This paper introduces a mesh-free operator learning framework for computing the fine-scale solution of multiscale PDEs. The proposed framework is trained using (i) coarse-scale solutions, which are inexpensive to obtain, and (ii) a limited number of observations of the fine-scale solution, to approximate the fine-scale solution at any desired location within the domain. Additionally, when the observations are noisy, we designed a Bayesian, multiscale operator learning approach that can reliably predict fine-scale solutions. Finally, we demonstrated the effectiveness and reliability of the proposed framework using a 1D elliptic equation, 2D elliptic equations with one fast variable and multiple scales, and the radiative transfer equation. The results confirmed that the proposed framework can work as a multiscale, mesh-free solver. In future work, we plan to design a DON that is invariant to input discretization. This will enable more effective capturing of derivatives by having patches with different discretizations.}

\section{Acknowledgement}
Z. Zhang was supported in part by AFOSR MURI FA9550-21-1-0084. H. Schaeffer was supported in part by AFOSR MURI FA9550-21-1-0084, NSF DMS-2208339, and an NSF CAREER Award DMS-2331100.

\bibliographystyle{abbrv}
\bibliography{references}
\end{document}